\documentclass[journal,transmag]{IEEEtran}
\usepackage{graphicx}
\usepackage{amsmath,amssymb,bm}
\usepackage{booktabs}
\usepackage{cite}   

\usepackage{tikz}
\usepackage{eso-pic}
   
\usetikzlibrary{arrows.meta}
\definecolor{coregray}{HTML}{E1E0D9} 
\definecolor{phaseA}{HTML}{2A78D6}
\definecolor{phaseB}{HTML}{EB6834}

\newcommand{\vH}{\bm{H}}
\newcommand{\vB}{\bm{B}}
\newcommand{\vM}{\bm{M}}
\newcommand{\vHr}{\bm{H}_{\mathrm{rev}}}
\newcommand{\vHp}{\bm{H}_{\mathrm{rev,p}}}
\newcommand{\Hr}{H_{\mathrm{rev}}}
\newcommand{\vr}{\bm{r}}

\newcommand{\vz}{\bm{z}}
\newcommand{\vR}{\bm{R}}
\newcommand{\vd}{\bm{d}}
\newcommand{\vF}{\bm{F}}
\newcommand{\vs}{\bm{s}}
\newcommand{\vy}{\bm{y}}
\newcommand{\vphi}{\bm{\phi}}
\newcommand{\hu}{\hat{\bm{h}}}
\newcommand{\hr}{\hat{\bm{r}}}
\newcommand{\tA}{\mathbf{A}}
\newcommand{\tP}{\mathbf{P}}
\newcommand{\tJ}{\mathbf{J}}
\newcommand{\tI}{\mathbf{I}}
\newcommand{\tK}{\mathbf{K}}
\newcommand{\tG}{\mathbf{G}}
\newcommand{\tS}{\mathbf{S}}
\newcommand{\tW}{\mathbf{W}}
\newcommand{\tC}{\mathbf{C}}
\newcommand{\tQ}{\mathbf{Q}}
\newcommand{\tY}{\mathbf{Y}}
\newcommand{\tE}{\mathbf{E}}
\newcommand{\tX}{\mathbf{X}}
\newcommand{\tZ}{\mathbf{0}}
\newcommand{\sym}{\operatorname{sym}}
\newcommand{\skw}{\operatorname{skew}}
\newcommand{\tr}{\operatorname{tr}}
\newcommand{\adj}{\operatorname{adj}}
\newcommand{\Man}{M_{\mathrm{an}}}
\newcommand{\wg}{w_{\mathrm g}}

\begin{document}
\AddToShipoutPictureFG*{\AtPageLowerLeft{\put(0,\LenToUnit{0.32in}){\makebox[\paperwidth]{\parbox{\textwidth}{\centering\footnotesize
This work has been submitted to the IEEE for possible publication.\\ Copyright may be transferred without notice, after which this version may no longer be accessible.}}}}}

\title{Loss of Monotonicity in the Vector Play Model and Consequences for Newton-Type Solvers}
\author{\IEEEauthorblockN{Dominik Cibicki\IEEEauthorrefmark{1}, Klaus Roppert\IEEEauthorrefmark{1}, and Manfred Kaltenbacher\IEEEauthorrefmark{1}}
\IEEEauthorblockA{\IEEEauthorrefmark{1}Institute of Fundamentals and Theory in Electrical Engineering, TU Graz, Graz, Austria}% <-this % stops an unwanted space
\thanks{Manuscript received October, 2026; Corresponding author: Dominik Cibicki (email: dominik.cibicki@tugraz.at).}}

\markboth{Journal of \LaTeX\ Class Files,~Vol.~14, No.~8, October 2026}%
{Cibicki \MakeLowercase{\textit{et al.}}: Loss of Monotonicity in the Vector Play Model and Consequences for Newton-Type Solvers}

\IEEEtitleabstractindextext{
\begin{abstract}
The vector play model (VPM) replaces the implicit minimization of the energy-based hysteresis
model (EBM) by an explicit update. Its tangent contains the product of the Hessians of two unrelated
potentials and is unsymmetric. Hence, the VPM has no co-energy density at frozen memory, and the
skew part makes $\vH\mapsto\vB$ lose monotonicity in a band just outside the pinning circle, into which
small time steps drive the sliding states. In composite models, the determinant of the tangent can change sign,
which explains published failures of the local inversion $\vB\mapsto\vH$, and where it is negative,
$\vH\mapsto\vB$ is not injective. In the scalar potential formulation,
the exact tangent keeps a descent guarantee, whereas symmetric surrogates of the tangent, such as the
symmetrized tangent or Broyden--Fletcher--Goldfarb--Shanno (BFGS) updates, lose it and become ascent
directions where the skew part dominates, as computations of TEAM Problem~32 confirm. It is therefore
recommended to use solvers with an unsymmetric tangent.

\end{abstract}
\begin{IEEEkeywords}
Energy-based hysteresis, finite element method, monotone operators, Newton's method,
vector play model.
\end{IEEEkeywords}}

\maketitle
\IEEEdisplaynontitleabstractindextext
\IEEEpeerreviewmaketitle

% =====================================================================

\section{Introduction}
\IEEEPARstart{V}{ector} hysteresis models that comply with thermodynamics, in particular the energy-based model
(EBM) of~\cite{franccois2013energy}, are increasingly used
in finite element (FE) analyses of electrical devices
\cite{Prigozhin2016,KaltenbacherRoppert2025,Domenig2024}. For the solver, the decisive property
of the EBM is variational: at frozen memory, its forward operator $\vH\mapsto\vB$, relating the magnetic field strength $\vH$ and the magnetic flux density $\vB$, is the gradient
of a convex co-energy density and strongly monotone~\cite{EggerTMAG2025,EggerEngertsberger2025}.
This guarantees unique solvability of the discrete scalar potential formulation and global,
mesh-independent convergence of damped Newton and local quasi-Newton methods
\cite{EggerEngertsberger2025,Egger2025CMA}. The price is a non-smooth convex minimization at every
quadrature point and for every pinning force~\cite{KaltenbacherRoppert2025,Egger2025eval}.

The vector play model (VPM) replaces this minimization by a closed-form geometric update
\cite{Prigozhin2016,Jacques2016,KaltenbacherRoppert2025}, which shortens the local inversion by up to
two orders of magnitude \cite{Jacques2018} and the simulation of TEAM Problem~32 by a factor
of about 2.4 \cite{Domenig2026}. However, convergence failures of the VPM in settings where the EBM converges have been reported \cite{Jacques2016,Jacques2018}, in particular for the material inversion $\vB \mapsto \vH$. They were attributed to the kinks
of the update at the pinning threshold and to near-singular susceptibilities in saturation
\cite{Jacques2016}.

We identify a structural cause of these failures that the kinks do not explain, for the isotropic VPM
driven by $\vH$ \cite{Bergqvist1997,Prigozhin2016} in 2-D and the reduced scalar potential formulation.
First, the VPM tangent contains the product of the Hessians of two unrelated potentials, and no
co-energy density exists (Sections~\ref{sec:tangent}, \ref{sec:asym}). Second, $\vH\mapsto\vB$ loses
monotonicity in a closed-form band near the pinning circle, into which time-step refinement drives the
sliding states (Sections~\ref{sec:mono}, \ref{sec:band}). Third, in composite models, the determinant of
the tangent can change sign, which explains the failures of the local inversion
in~\cite{Jacques2016,Jacques2018}, and $\vH\mapsto\vB$ is then not injective, i.e., the inversion has several solutions
(Sections~\ref{sec:comp}, \ref{sec:local}). Fourth, in the scalar potential formulation, the exact
tangent keeps a descent guarantee whereas symmetric surrogates lose it, as computations of TEAM
Problem~32 confirm (Section~\ref{sec:fe}).

% =====================================================================

\section{The Vector Play Model as a Material Law}
In each cell of the EBM and the VPM, the field splits into a reversible and an irreversible part,
$\vH=\vHr+\vH_{\mathrm{irr}}$. The magnetization $\vM=\Man(\Hr)\hu$, $\Hr:=\|\vHr\|$,
$\hu:=\vHr/\Hr$, follows from a concave anhysteretic curve $\Man\in C^1$ with $\Man(0)=0$,
$\Man'>0$, and $\vB=\mu_0(\vH+\vM)$ with vacuum permeability $\mu_0$. Throughout, $\|\cdot\|$ is the
Euclidean norm of vectors and $\|\cdot\|_2$ the spectral norm of matrices.
The two models differ only in the direction of $\vH_{\mathrm{irr}}$
\cite{Prigozhin2016,Jacques2016},
\begin{align}
\text{EBM:}&\quad\vH_{\mathrm{irr}}=\kappa\frac{\vM-\vM_{\mathrm p}}{\|\vM-\vM_{\mathrm p}\|},\nonumber\\
\text{VPM:}&\quad\vH_{\mathrm{irr}}=\kappa\frac{\vH-\vHp}{\|\vH-\vHp\|},\label{eq:dir}
\end{align}
where $\kappa$ is the pinning force, and $\vM_{\mathrm p}$ and $\vHp$ are the magnetization and the
reversible field of the previous load step.
The EBM direction thus depends on the unknown $\vM$, making it implicit, while
the VPM direction is known.

\subsection{Update Rule and Composite Model}\label{sec:model}
The explicit update rule of the VPM reads~\cite{Jacques2018,KaltenbacherRoppert2025}
\begin{align}
\vHr(\vH)&=\begin{cases}\vHp, & r\le\kappa\ \ \text{(pinned)},\\
\vH-\kappa\hr, & r>\kappa\ \ \text{(sliding)},\end{cases}\label{eq:update}\\
\vB(\vH)&=\mu_0\big(\vH+\Man(\Hr)\hu\big),\label{eq:B}
\end{align}
where $\vr:=\vH-\vHp$, $r:=\|\vr\|$, $\hr:=\vr/r$. The circle $r=\kappa$ around $\vHp$ is the pinning circle. In the composite model, necessary for modeling minor loops, the VPM superposes cells $k$ with weights $\omega_k>0$, $\sum_k\omega_k=1$, each with its
own $\kappa_k$ and $\vH_{\mathrm{rev,p},k}$. Hence
\begin{equation}
\vB=\mu_0\Big(\vH+\sum_{k}\omega_k\Man(H_{\mathrm{rev},k})\,\hu_k\Big).\label{eq:Bcomp}
\end{equation}
A cell with $\kappa_k=0$ is purely reversible. 
\subsection{Potentials and Tangent}\label{sec:tangent}
For a sliding cell, $\vHr$ and $\vM$ are gradients,
\begin{align}
\vHr&=\partial_{\vH}\Psi, & \Psi&:=\tfrac12\|\vH\|^2-\kappa r,\nonumber\\
\vM&=\partial_{\vHr}\wg, & \wg&:=\textstyle\int_0^{\Hr}\Man(s)\,\mathrm ds .\label{eq:pot}
\end{align}
In the sliding region, both potentials are locally convex~\cite{Prigozhin2016,MielkeRoubicek2015} but
unrelated: $\wg$ is constitutive, since $\mu_0(\frac12\|\vH\|^2+\wg(\vH))$ is the co-energy density of the anhysteretic
material, whereas $\Psi$ is geometric and only knows $\kappa$ and the memory. While the forward operator
of the EBM is the gradient of a single convex co-energy density $w^*$ such that
$\vB=\partial_{\vH}w^*$~\cite{EggerTMAG2025}, the VPM composes the gradients \eqref{eq:pot},
\begin{equation}
\vB(\vH)=\mu_0\big(\vH+(\partial_{\vHr}\wg\circ\partial_{\vH}\Psi)(\vH)\big),\label{eq:comp}
\end{equation}
which in general is not the gradient of a potential on its own (Section~\ref{sec:asym}). Let $\delta:=\kappa/r$, $\gamma:=1-\delta$, the secant and differential susceptibilities
$\alpha:=\Man(\Hr)/\Hr$, $\tau:=\Man'(\Hr)$, $\beta:=\tau-\alpha\le0$ and $\theta$ the signed angle from $\hr$ to $\hu$.
The Hessians are
\begin{align}
\tA&:=\partial^2_{\vHr}\wg=\alpha\tI+\beta\hu\otimes\hu,\nonumber\\
\tP&:=\partial^2_{\vH}\Psi=\gamma\tI+\delta\hr\otimes\hr .\label{eq:AP}
\end{align}
Both $\tA$ and $\tP$ are symmetric positive definite (SPD) for $r>\kappa$, with eigenvalues $\tau,\alpha$ and $1,\gamma$, and hence
$\det(\tA\tP)=\alpha\tau\gamma$. Pinned cells have $\tP_k=\tZ$, since $\vH_{\mathrm{rev},k}$ stays fixed, and
reversible ones $\tP_k=\tI$, since $\vH_{\mathrm{rev},k}=\vH$. On the
pinning circle, $\tP_k$ jumps from $\tZ$ to $\hr_k\otimes\hr_k$, which is the kink of the update. The
tangent is the weighted sum of the cell tangents $\tJ_k:=\mu_0(\tI+\tA_k\tP_k)$ \cite{Jacques2018},
\begin{equation}
\tJ:=\partial_{\vH}\vB=\sum_k\omega_k\tJ_k
=\mu_0\Big(\tI+\sum_k\omega_k\tA_k\tP_k\Big),\label{eq:J}
\end{equation}
with $J_{ij}=\partial B_i/\partial H_j$. To our knowledge, these factors have not been identified as Hessians of the
unrelated potentials \eqref{eq:pot} before.
In what follows, quantities of cell $k$ carry the index $k$, which we omit
for a single cell.

For the isotropic VPM in 3-D, the single-cell results below hold in the plane spanned by $\hr$ and
$\hu$, which $\tA$ and $\tP$ leave invariant. For anisotropic VPMs, we do not expect this reduction, and in
particular the closed forms below, to hold.

% =====================================================================

\subsection{Asymmetry and Non-Conservativity}\label{sec:asym}
Let $\sym\tY:=\frac12(\tY+\tY^T)$, $\skw\tY:=\frac12(\tY-\tY^T)$ and the commutator
$[\tY_1,\tY_2]:=\tY_1\tY_2-\tY_2\tY_1$ \cite{horn2012matrix}. Then $\vB$, $\tJ$, $\skw\tJ$ and $\sym\tJ$ are weighted sums
over the cells, whereas eigenvalues and determinants are not. Only sliding cells contribute to
$\skw\tJ$. For each of them, $(\tA\tP)^T=\tP\tA$ and thus $\skw(\tA\tP)=\frac12[\tA,\tP]$. Since the
multiples of $\tI$ in \eqref{eq:AP} commute, only $[\hu\otimes\hu,\hr\otimes\hr]=\cos\theta(\hu\otimes\hr-\hr\otimes\hu)
=\cos\theta\sin\theta\tE$ remains, with the 2-D rotation generator
$\tE:=\bigl(\begin{smallmatrix}0&-1\\1&0\end{smallmatrix}\bigr)$. Hence
\begin{equation}
\skw\tJ=\frac{\mu_0}{4}\sum_k\omega_k\beta_k\delta_k\sin2\theta_k\tE .\label{eq:skew}
\end{equation}
A sliding cell is symmetric if and only if $\beta=0$, $\kappa=0$ or $\sin2\theta=0$, i.e., for unidirectional
and stationary rotational loading, where the VPM agrees with the EBM \cite{Prigozhin2016}. In these
cases, $\hr$ is an eigenvector of $\tA$. Furthermore, each contribution is largest just beyond the pinning
circle, where $\delta_k\to1$.

\begin{figure*}[t]
\centering
\includegraphics[width=0.9\textwidth,trim=0 6.0 0 0,clip]{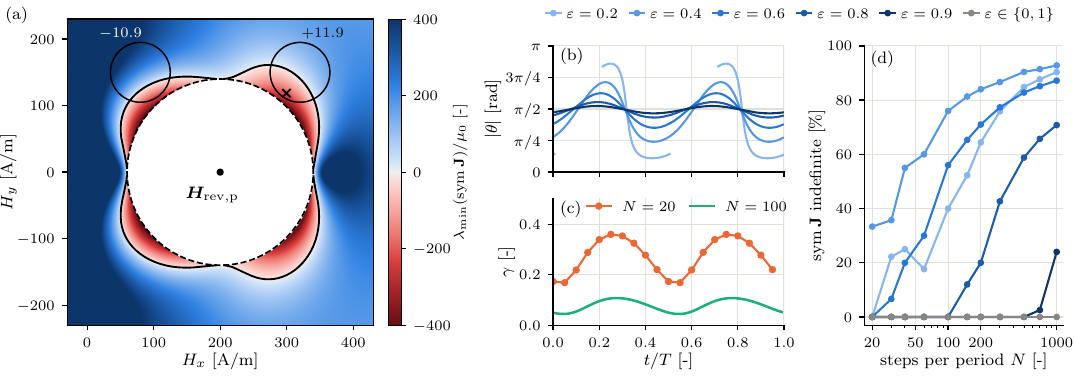}
\caption{Single cell, $\Man=(2M_{\mathrm s}/\pi)\arctan(\Hr/a)$, $M_{\mathrm s}=1.23$~MA/m,
$a=38$~A/m, $\kappa=140$~A/m \cite{Domenig2024}. (a)~$\lambda_{\min}(\sym\tJ)/\mu_0$, clipped to $\pm400$, for $\vHp=(200,0)^T$~A/m with blank pinned
interior. Thick line: $\lambda_{\min}=0$, dashed: pinning circle, thin circles: circulation \eqref{eq:circ} in J/m$^3$,
cross: $\vH=\vHp+155$~A/m$\,(\cos50^\circ,\sin50^\circ)^T$, $\lambda(\sym\tJ)/\mu_0=\{-201,3541\}$.
(b)--(d)~Settled period of
$\vH(t)=H_0(\cos2\pi t,\,\varepsilon\sin2\pi t)$, $H_0=300$~A/m: (b)~$|\theta|$ (gaps: pinned), (c)~$\gamma$ for $\varepsilon=0.6$ and
$N=20,100$, (d)~share of sliding steps with indefinite $\sym\tJ$ over~$N$.}
\label{fig:material}
\end{figure*}

The asymmetry has a global counterpart: On every disk $D$ of smooth states, i.e., with all
$\vH_{\mathrm{rev},k}\ne\bm 0$ and off all pinning circles, $\vB$ is continuously differentiable and the
set of sliding cells is fixed. Green's theorem with $J_{21}-J_{12}=2(\skw\tJ)_{21}$ then yields
\begin{equation}
\oint_{\partial D}\vB\cdot\mathrm d\vH
=\frac{\mu_0}{2}\iint_D\sum_k\omega_k\beta_k\delta_k\sin2\theta_k\,\mathrm dA .
\label{eq:circ}
\end{equation}
A co-energy density with $\partial_{\vH}w^*=\vB$ on such a disk would make every such circulation vanish,
or equivalently, $\tJ=\partial^2_{\vH}w^*$ would be symmetric \cite{marsden2012vector}. In Fig.~\ref{fig:material}(a), the circulations
of a single cell along the two thin circles of radius 45~A/m in the sliding region are $-10.9$ and
$11.9$~J/m$^3$ instead of zero. Hence, unlike the EBM \cite{EggerTMAG2025}, the VPM is not
conservative at frozen memory.

% =====================================================================

\subsection{Spectrum and Loss of Monotonicity}\label{sec:mono}
For a sliding cell, $\tA\tP=\tA^{1/2}\big(\tA^{1/2}\tP\tA^{1/2}\big)\tA^{-1/2}$ is similar to an
SPD matrix \cite{horn2012matrix} with spectral norm at most $\|\tA\|_2\|\tP\|_2=\alpha\le\Man'(0)$,
by concavity of $\Man$.
Therefore, its eigenvalues are real and positive, although $\tA\tP$ is unsymmetric, and every cell
tangent has real eigenvalues with
\begin{equation}
\mu_0\le\lambda_i(\tJ_k)\le\mu_0\big(1+\Man'(0)\big).\label{eq:spec}
\end{equation}
Hence, a single-cell tangent is always invertible with $\det\tJ\ge\mu_0^2$. Monotonicity of
$\vH\mapsto\vB$, however, is equivalent to $\sym\tJ\succeq0$ \cite{Zeidler1990}. Since $\sym\tJ$ is
the weighted sum of the $\sym\tJ_k$, the composite is monotone wherever all its cells are. For a
single cell, using the two identities $\det\tY=\det(\sym\tY)+\|{\skw\tY}\|_2^2$ and
$\det(\tI+\tY)=1+\tr\tY+\det\tY$, which hold for every $\tY\in\mathbb R^{2\times2}$ since the skew part
of a $2\times2$ matrix is a multiple of $\tE$, we obtain
\begin{align}
\det(\sym\tJ)&=\det\tJ-\|{\skw\tJ}\|_2^2,\nonumber\\
\frac{\det\tJ}{\mu_0^2}&=1+\tr(\tA\tP)+\det(\tA\tP).\label{eq:detsym}
\end{align}
Since $\tr(\sym\tJ)=\tr\tJ>0$, $\sym\tJ$ is indefinite if and only if $\det(\sym\tJ)<0$, i.e.,
$\|{\skw\tJ}\|_2^2>\det\tJ\ge\mu_0^2$. Hence, only the skew part can destroy monotonicity, and it
must outweigh the determinant of the full tangent. In one dimension, the VPM reduces to the monotone
scalar play operator \cite{KrejciMonteiro2025}, and the EBM is monotone in every dimension
\cite{EggerTMAG2025}. The VPM's loss of monotonicity is therefore a vectorial artifact of the update
direction $\hr$ in \eqref{eq:dir}. Along the eigenvector $\vz$ of a negative eigenvalue of
$\sym\tJ$, the incremental permeability is negative, such that a small increment $\Delta\vH\parallel\vz$ yields
$\Delta\vH\cdot\Delta\vB\approx\Delta\vH\cdot(\sym\tJ)\,\Delta\vH<0$, e.g., at the cross in
Fig.~\ref{fig:material}(a).

% =====================================================================

\subsection{Critical Band and Time Step}\label{sec:band}
For a single cell, we insert $\det(\tA\tP)=\alpha\tau\gamma$ into \eqref{eq:detsym} and exploit
\eqref{eq:skew} to obtain $\|{\skw(\tA\tP)}\|_2^2=\beta^2\delta^2\sin^22\theta/16$. With
$q:=\alpha/\tau\ge1$, we have $\beta^2=(\alpha-\tau)^2=16\alpha\tau\rho(q)$ and thus
\begin{equation}
\frac{\det(\sym\tJ)}{\mu_0^2}=1+\tr(\tA\tP)+\alpha\tau\big(\gamma-\rho(q)\delta^2\sin^22\theta\big),\label{eq:detband}
\end{equation}
with $\rho(q):=(q-1)^2/(16q)$. The bracket compares the stabilizing $\det(\tA\tP)\propto\gamma$ with the destabilizing skew part.
Since $\tr(\sym\tJ)>0$, $\vH\mapsto\vB$ is locally monotone exactly if \eqref{eq:detband} is
positive, i.e., if the bracket exceeds $-(1+\tr(\tA\tP))/(\alpha\tau)$. This
bound is negative and of the order $1/\tau$, where $\tau=\Man'(\Hr)$ is the relative differential
susceptibility (typically $10^2$ to $10^4$ in iron). Neglecting it and using $\delta=1-\gamma$ gives
\begin{equation}
\gamma>\rho(q)(1-\gamma)^2\sin^22\theta,\label{eq:band}
\end{equation}
which is a sufficient condition for monotonicity and the exact condition for $\sym(\tA\tP)\succeq0$,
i.e., without the vacuum part $\mu_0\tI$. It separates saturation ($q$), loading geometry ($\theta$)
and time step ($\gamma$). Since the left-hand side grows with $\gamma=1-\kappa/r$, \eqref{eq:band}
fails in a band $0<\gamma<\gamma_{\mathrm{crit}}(q,\theta)$ just outside the pinning circle whenever
$\sin2\theta\ne0$. Up to the neglected term of order $1/\tau$, this is the band of lost
monotonicity. For example, $\gamma_{\mathrm{crit}}=0.217$, i.e., $r/\kappa<1.28$, at $\Hr=200$~A/m
and $\theta=45^\circ$ (see Fig.~\ref{fig:material}(a)).

The time step enters through $\gamma$. After every update, $\|\vH-\vHr\|\le\kappa$ holds. Hence, in the
next step with increment $\Delta\vH_n$, $r_n\le\kappa+\|\Delta\vH_n\|$ holds, and therefore
\begin{equation}
\gamma_n=1-\kappa/r_n\le\|\Delta\vH_n\|/\kappa=O(\Delta t).\label{eq:gammadt}
\end{equation}
In contrast, $q$ and $\theta$ converge to values set by the load. Therefore, time-step refinement drives
every sliding state into its band, and monotonicity is lost unless $\sin2\theta$ is small.
Fig.~\ref{fig:material}(b)--(d) shows this for elliptical loading with axis ratio $\varepsilon$. For
$\varepsilon=0.6$, the share of sliding steps with indefinite $\sym\tJ$ grows from 0 to 87~\% between
$N=20$ and 1000 steps per period, and only unidirectional ($\varepsilon=0$) and circular
($\varepsilon=1$) loading remain monotone. In the composite model, every sliding cell has its own band,
and monotonicity is lost only if the negative part of a cell, weighted by $\omega_k$, outweighs
$\mu_0\tI$ and the other cells, in particular the reversible one.

% =====================================================================

\subsection{Determinant of the Composite Model}\label{sec:comp}
The weighted sum \eqref{eq:J} preserves monotonicity, but not the spectral bound \eqref{eq:spec}.
Although the cell tangents have positive eigenvalues, they are unsymmetric, and their weighted sum can be
singular. Let $\tX_k:=\omega_k\tA_k\tP_k$, such that $\tJ/\mu_0=\tI+\sum_k\tX_k$, and let
$\adj\tX:=(\tr\tX)\tI-\tX$ be the adjugate. Repeated use of
$\det(\tY_1+\tY_2)=\det\tY_1+\det\tY_2+\tr(\adj\tY_1\,\tY_2)$ for $\tY_1,\tY_2\in\mathbb R^{2\times2}$ \cite{horn2012matrix}
generalizes \eqref{eq:detsym} to
\begin{equation}
\frac{\det\tJ}{\mu_0^2}=1+\sum_k\big(\tr\tX_k+\det\tX_k\big)+\sum_{k<l}\tr(\adj\tX_k\,\tX_l).
\label{eq:detcomp}
\end{equation}
Since every $\tX_k$ has real non-negative eigenvalues, the single-cell terms are non-negative, with
$\det\tX_k=\omega_k^2\alpha_k\tau_k\gamma_k$. For symmetric positive semidefinite $\tX_k$, as in the
EBM, the mixed terms are non-negative as well, and $\det\tJ\ge\mu_0^2$ follows. 

For the VPM, however, only
the mixed terms of cells that are not pinned can be negative, and only if one of them is in its band
or if their skew parts have opposite sense. Since $\det\tJ=\det(\sym\tJ)+\|{\skw\tJ}\|_2^2$,
$\det\tJ<0$ requires $\sym\tJ$ to be indefinite. Hence, sign changes of $\det\tJ$ lie in the
non-monotone region and need at least two cells that are not pinned. The determinant changes sign
either continuously, through $\det\tJ=0$ inside the sliding region, or by a jump across a pinning
circle, where the mixed terms of cell $k$ appear at once.\looseness=-1

% =====================================================================

\subsection{Material-Point Test: Local Inversion}\label{sec:local}
Vector potential formulations require the inversion $\vB\mapsto\vH$. For the EBM,
convex duality yields it as a minimization \cite{EggerTMAG2025}, whereas for the VPM one solves
\cite{Moll2025,Domenig2026}
\begin{equation}
\vF(\vH):=\vB(\vH)-\vB^*=\bm 0 .\label{eq:res}
\end{equation}
The Newton--Raphson (NR) method iterates $\vH_{l+1}=\vH_l+\eta_l\vd_l$ with $\vd_l:=-\tJ^{-1}\vF$, either
undamped ($\eta_l=1$) or with backtracking, i.e., $\eta_l=2^{-m}$ with the smallest $m$ that decreases
$\|\vF\|$. For $f:=\frac12\|\vF\|^2$, we may write
\begin{equation}
\partial_{\vH}f\cdot\vd_l=\vF\cdot\tJ\vd_l=-\|\vF\|^2<0\label{eq:descentloc}
\end{equation}
wherever $\tJ$ is invertible \cite{DennisSchnabel1996}, regardless of its symmetry or the sign of
$\det\tJ$. Quasi-Newton (QN) methods replace $\tJ$ by matrices $\tQ_l$ that satisfy the secant
condition $\tQ_{l+1}\vs_l=\vy_l$ with $\vs_l:=\vH_{l+1}-\vH_l$ and $\vy_l:=\vF(\vH_{l+1})-\vF(\vH_l)$
\cite{DennisSchnabel1996,NocedalWright2006}. Broyden's update can represent $\skw\tJ$ but has no
descent guarantee. The Broyden--Fletcher--Goldfarb--Shanno (BFGS) update, used in \cite{Jacques2018,Egger2025CMA}, is symmetric and stays positive definite if and only if
\begin{equation}
\vs_l\cdot\vy_l=\int_0^1\vs_l\cdot\sym\tJ(\vH_l+t\vs_l)\vs_l\,\mathrm dt>0 ,\label{eq:curvloc}
\end{equation}
which is the secant form of monotonicity.

\begin{table}[t]
\caption{Local material inversion \cite{Jacques2018}: number of excitations with at least one failed
step}
\label{tab:jacques}
\centering\footnotesize
\setlength{\tabcolsep}{3.4pt}
\begin{tabular}{@{}lccccc@{}}
\toprule
Trajectory class & \multicolumn{3}{c}{with backtracking} & \multicolumn{2}{c}{undamped}\\
\cmidrule(lr){2-4}\cmidrule(l){5-6}
(excitations) & NR & BFGS & Broyden & NR & BFGS\\
\midrule
\multicolumn{6}{@{}l}{\emph{VPM, 3 cells}}\\
\quad$\det\tJ<0$ (3) & 3 & 3  & 3  & 3  & 3\\
\quad$\det\tJ\ge0$, $\sym\tJ$ indef. (16) & 2 & 16 & 2  & 16 & 12\\
\quad$\sym\tJ\succ0$ (21) & 0 & 4  & 0  & 10 & 0\\
\quad all (40) & 5 & 23 & 5  & 29 & 15\\
\quad all (40), as in \cite{Jacques2018} & 5 & 23 & -- & 30 & 17\\
\multicolumn{6}{@{}l}{\emph{VPM, 1 cell}}\\
\quad$\det\tJ\ge0$, $\sym\tJ$ indef. (31) & 0 & 31 & 4  & 31 & 31\\
\quad$\sym\tJ\succ0$ (9) & 0 & 1  & 0  & 9  & 9\\
\multicolumn{6}{@{}l}{\emph{EBM} \cite{Jacques2018}}\\
\quad all (40) & 0 & 0 & -- & 30 & 0\\
\bottomrule
\multicolumn{6}{@{}p{0.97\columnwidth}@{}}{\scriptsize\raggedright Classes by the true trajectory: $\det\tJ<0$ at
some step, $\det\tJ\ge0$ throughout but $\sym\tJ$ indefinite (indef.) at some step, or $\sym\tJ\succ0$ at all steps.
A step fails if it does not converge within 200 iterations.
Row \cite{Jacques2018}: published values, with relaxation instead of backtracking. EBM: variational and
differential form, identical.}
\end{tabular}
\end{table}

\begin{figure}[t]
\centering
\includegraphics[width=1.0\columnwidth,trim=0 4.1 0 0,clip]{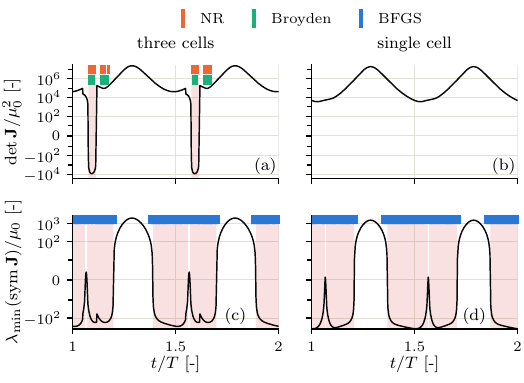}
\caption{Local inversion, $H_{\max}=200$~A/m, $\Phi=20^\circ$, second period: (a),~(b)~$\det\tJ/\mu_0^2$,
(c),~(d)~$\lambda_{\min}(\sym\tJ)/\mu_0$ on the true trajectory, three cells and one cell, negative values
shaded. Ticks mark failed steps with backtracking.}
\label{fig:inversion}
\end{figure}

We repeat the inversion tests of \cite{Jacques2018} with three cells, $(\omega_k,\kappa_k)=(0.1,0)$,
$(0.65,16\,\mathrm{A/m})$, $(0.25,47\,\mathrm{A/m})$ and $\mu_0\Man=1.22\,\mathrm T\cdot\tanh(\Hr/65\,\mathrm{A/m})$,
and add a single cell with $\kappa=16$~A/m and Broyden's method (Table~\ref{tab:jacques}). The 40
excitations $\vH(t)=H_{\max}(\cos2\pi t,\cos(2\pi t-\Phi))^T$, $0\le t\le2$, with $N=500$ steps per
period, are all $4\times10$ combinations of $H_{\max}\in\{50,100,150,200\}$~A/m and
$\Phi\in\{0^\circ,10^\circ,\dots,90^\circ\}$, from unidirectional to circular loading. The (quasi-)Newton iterations stop at
$\|\vF(\vH_l)\|\le10^{-7}$~T or after 200 iterations. After a failed step, the inversion
continues from the true state in our implementation.

Undamped NR fails for all models, since it cycles at the kinks \cite{Jacques2016}. With backtracking
or BFGS, only the VPM fails, which was attributed to the asymmetry of $\tJ$ \cite{Jacques2018}. By
\eqref{eq:descentloc}, however, damped NR can only stall where $\tJ$ is (nearly) singular or at a kink.
Kinks alone are harmless, as the single cell shows, but in composite models $\det\tJ$ changes sign
(Section~\ref{sec:comp}). Accordingly, damped NR fails exactly in the three excitations with
$\det\tJ<0$ on the trajectory and in the two adjacent ones ($H_{\max}=200$~A/m, $\Phi=10^\circ$ to
$50^\circ$), which are the five of \cite{Jacques2018}, and 168 of the 173 failed steps reach
$\det\tJ\le0$. There, the negative determinant stems from the mixed term of the reversible cell and
the cell with $\kappa=16$~A/m in its band.\looseness=-1

Moreover, at each of the 126 steps with $\det\tJ<0$, a
Newton search from a dense grid of starting points finds three solutions of \eqref{eq:res}, up to
21~A/m apart, implying non-injectivity of $\vH\mapsto\vB$. The damped NR, started from the previous step, converges at 78 of these steps to a wrong
solution, which the test on $\|\vF\|$ cannot detect and Table~\ref{tab:jacques} counts as converged. For
the single cell, $\det\tJ\ge\mu_0^2$ (Fig.~\ref{fig:inversion}(b)), and damped NR never fails despite
the same kinks. BFGS fails mainly where $\sym\tJ$ is indefinite (Fig.~\ref{fig:inversion}(c),~(d)),
since \eqref{eq:curvloc} can fail there, also for a single cell, whereas Broyden's update fails for
three cells exactly where NR does. For the EBM, strong monotonicity \cite{EggerTMAG2025} rules out both mechanisms.

% =====================================================================

\section{Finite Element Context: TEAM Problem 32}\label{sec:fe}
\subsection{Field Problem and Solvers}\label{sec:fesolver}
\begin{figure}[t]
\centering
\begin{minipage}[b]{0.555\columnwidth}\centering
% v36 (2026-10-25, Dominik): scale 0.20898 mm per mm = printed scale of the core in Fig. 4 (fig_t32_theta_v36.pdf,
%   printed by fig_t32_theta_v36.py), so that both cores have exactly the same size (v28-v33: 0.192 mm, which
%   matched the Fig. 4 of v28 but not the one with the gap between the cores of v30); core fill coregray is
%   the palette neutral #E1E0D9 (defined in the .tex).
% TEAM Problem 32 -- compact version for the one-row Fig. 3 (v23, 2026-10-22; v28: scale 0.192 mm = core of Fig. 4;
% v33: label 120 rotated along its dimension line, it covered the left winding): no text labels,
% limb width 30 mm and window width 42.25 mm as in the FE mesh (geo_mesh2D_t32_ns*.mat), windings in
% phaseA/phaseB as in fig_t32_geometry_noprobes.tikz. Coordinates in mm, origin = lower left core corner.
\begin{tikzpicture}[
    x=0.20898mm, y=0.20898mm,
    >={Latex[length=0.9mm,width=0.65mm]},
    core/.style    = {line width=0.4pt, draw=black, fill=coregray},
    win/.style     = {line width=0.35pt, draw=black, fill=white},
    dimline/.style = {line width=0.25pt, draw=black, <->},
    extline/.style = {line width=0.25pt, draw=black},
    dimlbl/.style  = {font=\scriptsize, inner sep=1pt}
  ]
\def\Wc{174.5}\def\Hc{175}\def\tl{30}\def\ww{42.25}\def\yk{27.5}\def\hw{120}
\pgfmathsetmacro{\xa}{\tl}\pgfmathsetmacro{\xb}{\tl+\ww}
\pgfmathsetmacro{\xc}{\Wc-\tl-\ww}\pgfmathsetmacro{\xd}{\Wc-\tl}
\pgfmathsetmacro{\yb}{\Hc-\yk-\hw}\pgfmathsetmacro{\yt}{\Hc-\yk}
\pgfmathsetmacro{\mM}{\Wc/2}\pgfmathsetmacro{\mW}{\tl+\ww/2}
\def\bw{6}\def\bb{37.5}\def\bt{137.5}
\filldraw[core] (0,0) rectangle (\Wc,\Hc);
\filldraw[win]  (\xa,\yb) rectangle (\xb,\yt);
\filldraw[win]  (\xc,\yb) rectangle (\xd,\yt);
\fill[phaseA] (-10,\bb) rectangle ({-10+\bw},\bt);
\fill[phaseA] ({\xa+4},\bb) rectangle ({\xa+4+\bw},\bt);
\fill[phaseB] ({\xd-4-\bw},\bb) rectangle ({\xd-4},\bt);
\fill[phaseB] ({\Wc+10-\bw},\bb) rectangle ({\Wc+10},\bt);
% dimensions
\draw[extline] (0,\Hc) -- (0,200); \draw[extline] (\Wc,\Hc) -- (\Wc,200);
\draw[dimline] (0,195) -- (\Wc,195); \node[dimlbl,above] at (\mM,195) {174.5};
\draw[extline] (0,0) -- (-22,0); \draw[extline] (0,\Hc) -- (-22,\Hc);
\draw[dimline] (-18,0) -- (-18,\Hc); \node[dimlbl,left] at (-18,{\Hc/2}) {175};
\draw[dimline] (0,{\yb+12}) -- (\xa,{\yb+12}); \node[dimlbl,above] at ({\tl/2},{\yb+12}) {30};
\draw[dimline] (\mW,\yb) -- (\mW,\yt); \node[dimlbl,fill=white,rotate=90] at (\mW,{\Hc/2}) {120};  % v33: rotated, so that the label no longer covers the winding
\node[font=\footnotesize] at (\mM,{\yb/2}) {$\Omega_{\mathrm{Fe}}$};
\end{tikzpicture}
\end{minipage}\hfill
\begin{minipage}[b]{0.435\columnwidth}\centering
\includegraphics[width=\linewidth]{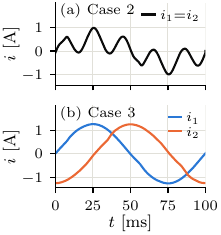}
\end{minipage}
\caption{TEAM Problem~32 \cite{Bottauscio2002}: iron core $\Omega_{\mathrm{Fe}}$ (in mm) with windings 1 (blue) and 2
(orange) inside the air box $\Omega$ of $245\times245$~mm (not depicted), and measured winding currents over one
period, (a)~Case~2, (b)~Case~3.}
\label{fig:t32geo}
\end{figure}

TEAM Problem~32 \cite{Bottauscio2002}, a benchmark for vector hysteresis models, is a three-limb laminated core with two windings of 90 turns, driven by
measured currents with period $T=100$~ms (Fig.~\ref{fig:t32geo}). In Case~2, both windings
carry the same current with a strong fifth harmonic, which leads to alternating flux. In Case~3, the
currents are almost sinusoidal and shifted by about $90^\circ$, which leads to rotational flux in the
T-joints. 

At every time step $t_n$, with the memory frozen at $t_{n-1}$,
the field in the 2-D cross-section $\Omega$ satisfies $\nabla\times\vH=\bm J_{\mathrm s}$ and
$\nabla\cdot\vB=0$, with source current density $\bm J_{\mathrm s}$, $\vB\cdot\bm n=0$ on
$\partial\Omega$, and the material law according to \eqref{eq:Bcomp} in the core $\Omega_{\mathrm{Fe}}$ and
$\vB=\mu_0\vH$ elsewhere.
Since the VPM is driven by $\vH$, we use the reduced scalar potential formulation \cite{Domenig2024}:
with $\nabla\times\vH_{\mathrm s}=\bm J_{\mathrm s}$, the ansatz $\vH=\vH_{\mathrm s}-\nabla\phi$ fulfills
Amp\`ere's law. The weak formulation therefore reads: find $\phi\in H^1(\Omega)/\mathbb R$ such that
\begin{equation}
\int_\Omega\vB(\vH_{\mathrm s}-\nabla\phi)\cdot\nabla v\,\mathrm d\Omega=0\qquad\forall v\in H^1(\Omega).
\label{eq:weak}
\end{equation}
For the EBM, \eqref{eq:weak} is the optimality condition of the convex functional
$\mathcal E(\phi)=\int_\Omega w^*(\vH_{\mathrm s}-\nabla\phi)\,\mathrm d\Omega$
\cite{EggerEngertsberger2025,Egger2025CMA}, whereas the VPM is variationally inconsistent: by \eqref{eq:circ}, no such functional exists.
With piecewise linear $\phi_h=\sum_i\phi_iv_i$, $\vphi:=(\phi_i)_i$, on triangles $\Omega_e$ and $\vH_{\mathrm s}$ at the
centroids, $\vH_e$, $\vB_e$ and $\tJ_e$ are elementwise constant, and the residual
$R_i(\vphi):=\int_\Omega\vB(\vH_{\mathrm s}-\nabla\phi_h)\cdot\nabla v_i\,\mathrm d\Omega$ has the tangent matrix
\begin{equation}
\tK:=-\partial_{\vphi}\vR=\sum_e|\Omega_e|\tG_e^T\tJ_e\tG_e,\label{eq:K}
\end{equation}
with the element gradient matrices $\tG_e$ \cite{Kaltenbacher2015}. Newton's method solves $\tK\Delta\vphi=\vR$ and takes the largest step $\eta\Delta\vphi$,
$\eta=2^{-m}$, $m\le10$, with $\|\vR(\vphi+\eta\Delta\vphi)\|\le(1-\sigma\eta)\|\vR(\vphi)\|$ (Armijo condition with $\sigma = 10^{-4}$
\cite{DennisSchnabel1996}), or else the one with the smallest residual. After convergence, the memory is updated.

Since assembly is linear, $\sym\tK$ and $\skw\tK$ follow from \eqref{eq:K} with $\sym\tJ_e$ and
$\skw\tJ_e$, and only elements in the core with sliding cells contribute to $\skw\tK$. With $\vF\to-\vR$ and
$\tJ\to\tK$, Section~\ref{sec:local} carries over, and the exact direction satisfies
$\partial_{\vphi} f\cdot\Delta\vphi=-\|\vR\|^2<0$. Local QN methods \cite{Egger2025CMA} replace each $\tJ_e$ by a
$2\times2$ matrix updated with $\Delta\vH_e$ and $\Delta\vB_e$. Hence, BFGS requires
$\Delta\vH_e\cdot\Delta\vB_e>0$, whereas Broyden keeps the skew part.

A symmetric surrogate, such as the symmetrized tangent, discards $\tW:=\skw\tK$. For
$\tS:=\sym\tK$, $\Delta\vphi_S:=\tS^{-1}\vR$ and $\hat{\vR}:=\vR/\|\vR\|$, we obtain
\begin{align}
-\partial_{\vphi} f\cdot\Delta\vphi_S/\|\vR\|^2&=1+\hat{\vR}\cdot\tC\hat{\vR},\nonumber\\
\tC&:=\tfrac12\big(\tW\tS^{-1}-\tS^{-1}\tW\big),\label{eq:comm}
\end{align}
where $1$ is the value of the exact direction and the second term is the price of discarding $\tW$.
Since $\tC$ is traceless with $\|\tC\|_2\le\|\tW\|_2/\lambda_{\min}(\tS)$ for $\tS\succ0$, the
surrogate can be an ascent direction where $\skw\tK$ dominates, even if $\sym\tK\succ0$. For any other
symmetric approximation $\tQ$ of $\tK$, e.g., assembled from BFGS, Davidon--Fletcher--Powell or symmetric
rank-one updates of $\tJ_e$ \cite{NocedalWright2006}, \eqref{eq:comm} holds with $\tQ$ in place of
$\tS$, except that the $1$ becomes $\hat{\vR}\cdot\tS\tQ^{-1}\hat{\vR}$. Hence, none of them has a descent guarantee where $\skw\tK$ dominates. Of these
updates, we computed only BFGS. By
\eqref{eq:skew} and \eqref{eq:gammadt}, $\|{\skw\tJ_e}\|_2$ is largest for $\sin2\theta$ far from zero
(rotational flux) and small $\gamma$ (small time steps), while by \eqref{eq:band} the contribution of
the sliding cells to $\lambda_{\min}(\sym\tJ_e)$ shrinks and eventually turns negative. Therefore,
symmetric surrogates deteriorate under time-step refinement, whereas \eqref{eq:descentloc} is
unaffected.

\begin{figure}[t]
\centering
\includegraphics[width=\columnwidth,trim=0 5.4 0 0,clip]{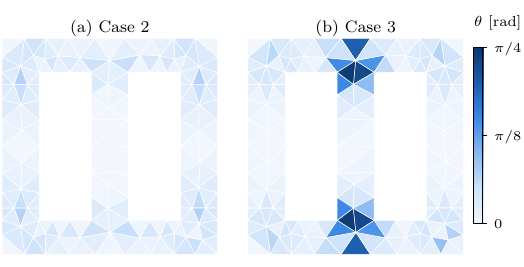}
\caption{TEAM Problem~32, mesh of Table~\ref{tab:dt} (170 elements in the core), 1000 steps per period: mean of the angle $\theta$
between $\hr$ and $\pm\hu$ over the sliding cells and the second period, (a)~Case~2, (b)~Case~3.}
\label{fig:t32mech}
\end{figure}

\subsection{Results}\label{sec:t32}
We compare four choices of the element tangents $\tJ_e$ in \eqref{eq:K}: the exact tangent, element-local
Broyden updates, the symmetrized tangent and element-local BFGS updates.
The core carries the 51 pinning forces (2.9~\% of the weight at $\kappa=0$) and the double-Langevin
anhysteretic curve identified in \cite{Domenig2024} from uniaxial measurements \cite{Domenig2026},
for which the VPM and the EBM coincide \cite{Prigozhin2016}. Every run starts from the demagnetized state and covers two periods. Newton stops
at $\|\vR\|\le10^{-8}\|\vR_0\|$, with $\vR_0$ the residual of the incoming iterate. A step that does not
converge within 300 iterations counts as failed, and the symmetrized tangent and BFGS stop at their first
failed step. Table~\ref{tab:dt} varies the time step on the mesh with 170 elements in the core, and
Table~\ref{tab:mesh} the mesh at 1000 steps per period. Both give the iterations averaged over the last
1.5 periods or, in parentheses, the time of the first step that stagnates above $10^{-6}\|\vR_0\|$, at
which the symmetrized tangent and BFGS stop.

\begin{table}[t]
\caption{TEAM Problem~32, two periods: Newton iterations or, in parentheses, time of the first failure, and
$\lambda_e:=\lambda_{\min}(\sym\tJ_e)$}
\label{tab:dt}
\centering\footnotesize
\setlength{\tabcolsep}{2.1pt}
\begin{tabular}{@{}lcccccc@{}}
\toprule
& \multicolumn{4}{c}{Case 3} & Case 3 & Case 2\\
& \multicolumn{4}{c}{identified} & uniform & identified\\
\cmidrule(lr){2-5}\cmidrule(lr){6-6}\cmidrule(l){7-7}
steps per period & 200 & 500 & 1000 & 2000 & 1000 & 1000\\
\midrule
exact & 6.3 & 5.8 & 5.5 & 5.3 & 5.0 & 5.3\\
Broyden & 9.1 & 7.1 & 6.4 & 6.2 & 5.3 & 7.9\\
symmetrized & 21.5 & (10.4\,ms) & (8.9\,ms) & (9.0\,ms) & 14.9 & 7.4\\
BFGS & 21.6 & (11.4\,ms) & (9.3\,ms) & (9.1\,ms) & 13.4 & 7.8\\
$\min_e\lambda_e/\mu_0$ & $-62$ & $-132$ & $-167$ & $-187$ & $27$ & $10.7$\\
\bottomrule
\end{tabular}
\end{table}

In Case~3, the exact tangent converges in every step with at most 10 iterations, and Broyden in all but
one with at most 40. The symmetrized tangent and BFGS need up to 290 iterations in converged steps. At 200
steps per period, they need 21.5 and 21.6 iterations on average instead of 6.3, and from 500 steps per
period on, they fail within 12~ms, stagnating above $10^{-3}$ and $10^{-5}\|\vR_0\|$ with ascent directions
in 91 to 98 and 16 to 51~\% of the 300 iterations. Although single element tangents are indefinite
($\min_e\lambda_e<0$), $\sym\tK$ stays positive definite in all runs, with smallest eigenvalue 1 relative
to the stiffness matrix of vacuum. Hence, the failures stem from $\tC$ in \eqref{eq:comm}, not from an
indefinite $\sym\tK$.

The last two columns of Table~\ref{tab:dt} are controls. In Case~2, where the local field alternates,
$\sym\tJ_e$ stays positive definite, and all four methods need 5 to 8 iterations. In Case~3 with the uniform
weights of \cite{Domenig2026}, $\kappa_k=0,35,\dots,140$~A/m and $\omega_k=0.2\ \forall k$, and the same anhysteretic curve,
$\sym\tJ_e$ stays positive definite as well ($\min_e\lambda_e=27\mu_0$), and no method fails at any of the
four time steps, although the symmetrized tangent and BFGS need 14.9 and 13.4 instead of 5.0 iterations at
1000 steps per period. Hence, rotational flux together with the identified distribution is critical.
The loss of monotonicity is also visible in the solver: \eqref{eq:curvloc} is violated in 125 element updates of BFGS
at 200 steps per period, but in none in Case~2 and with the uniform weights.

Fig.~\ref{fig:t32mech} shows that $\theta$ stays below 0.22~rad in Case~2, where the local field
alternates, but reaches 0.76~rad in the T-joints of Case~3, where it rotates. There, $\sin2\theta$, the
geometric factor of $\skw\tJ_e$ in \eqref{eq:skew}, is close to its maximum at $\theta=\pi/4$. The median
of $\gamma$ in Case~3 falls from 0.056 at 200 to 0.006 at 2000 steps per period, i.e., in proportion to
the time step.

\begin{table}[t]
\caption{TEAM Problem~32, Case~3, two periods, 1000 steps per period: mesh refinement, entries as in
Table~\ref{tab:dt}}
\label{tab:mesh}
\centering\footnotesize
\begin{tabular*}{\columnwidth}{@{\extracolsep{\fill}}lcccc@{}}
\toprule
elements in the core & 170 & 680 & 2720 & 10\,880\\
\midrule
exact & 5.5 & 5.7 & 5.9 & 6.4\\
Broyden & 6.4 & 6.6 & 7.4 & 9.3\\
symmetrized & (8.9\,ms) & 36.2 & 44.6 & (38.4\,ms)\\
BFGS & (9.3\,ms) & (11.0\,ms) & (12.1\,ms) & (11.7\,ms)\\
$\min_e\lambda_e/\mu_0$ & $-167$ & $-134$ & $-195$ & $-200$\\
\bottomrule
\end{tabular*}
\par\vspace{1pt}\parbox{\columnwidth}{\scriptsize\raggedright 170 elements: as in Table~\ref{tab:dt}.}
\end{table}

Table~\ref{tab:mesh} refines the mesh at 1000 steps per period. In Case~3, BFGS fails on every mesh within
14~ms, and the symmetrized tangent converges on the two intermediate meshes, with 36 and 45 iterations, but
fails at 38.4~ms on the finest one. The exact tangent and Broyden converge on all meshes, but their
iterations grow from 5.5 to 6.4 and from 6.4 to 9.3. Hence, in Case~3, convergence deteriorates under mesh
refinement, more than in Case~2, where the exact tangent needs 5.3 to 5.7 iterations.

% =====================================================================

\section{Conclusion}
The VPM composes the gradients of two unrelated potentials and is therefore variationally inconsistent, i.e., it has no potential of its own. Its
tangent is unsymmetric, and $\vH\mapsto\vB$ loses monotonicity in a band near the pinning circle that
small time steps reach. In composite models, $\det\tJ$ can change sign.
There, damped Newton stalls, which explains the failures of the local inversion, and where
$\det\tJ<0$, $\vH\mapsto\vB$ is not injective, the inversion has three solutions, and damped Newton often converges to a wrong solution. For Newton-type
solvers in the scalar potential formulation, we therefore recommend the exact unsymmetric tangent, or
Broyden's update, which needs no derivative of the material law but has no descent guarantee. Under rotational flux, the tangent should not be symmetrized, neither directly nor
by symmetric quasi-Newton updates such as BFGS, and smaller time steps do not help. Unlike for the EBM, mesh-independent convergence is not
guaranteed, and under rotational flux, convergence deteriorates under mesh refinement.

Additional considerations, such as the thermodynamic consistency of the VPM, i.e., whether its dissipation rate is non-negative and
maximal and whether it fits the variational framework of rate-independent systems \cite{MielkeRoubicek2015}, are left for future work.\looseness=-1

\section*{Acknowledgment}
This work was supported by the joint DFG/FWF Collaborative
Research Centre CREATOR (CRC – TRR361/F90) at TU
Darmstadt, TU Graz and JKU Linz.

\bibliographystyle{IEEEtran} 
\bibliography{vpm_refs} 

\end{document}